# Sifting Function Partition for the Goldbach Problem

Fu-Gao Song

College of Electronic Science & Technology, Shenzhen University, China

Shenzhen Key Laboratory of Micro-nano Photonic Information Technology, China

**Abstract.** All sieve methods for the Goldbach problem sift out all the composite numbers; even though, strictly speaking, it is not necessary to do so and which is, in general, very difficult. Some new methods introduced in this paper show that the Goldbach problem can be solved under sifting out only some composite numbers. In fact, in order to prove the Goldbach conjecture, it is only necessary to show that there are prime numbers left in the residual integers after the initial sifting! This idea can be implemented by using one of the three methods called "*sifting function partition by integer sort*", "*sifting function partition by intervals*" and "*comparative sieve method*", respectively. These are feasible methods for solving both the Goldbach problem and the problem of twin primes. An added bonus of the above methods is the elimination of the indeterminacy of the sifting functions brought about by their upper and lower bounds.

# 1. Introduction

In this paper, the term "*Proposition {1, b}*" means that *every large even number*

*can be expressed as the sum of a prime and the product of at most b primes.*

Rényi [1] proved in 1948, using Линник's large sieve [2], that there must be a positive constant $\eta_0$ such that the following estimation of the remainder term holds for any positive number $\eta < \eta_0$ and any $A > 0$:

$$\mathscr{R}(x,\eta) = \sum_{d \leq x^\eta} \mu^2(d) \max_{y \leq x} \max_{(l,d)=1} \left| \psi(y;l,d) - \frac{y}{\varphi(d)} \right| = O\left(\frac{x}{\ln^A x}\right); \qquad (1.1)$$

therefore proving Proposition {1, b} by Brun's sieve [3]. However, the values of $\eta_0$ and $b > 0$ were not determined. In fact, $\eta_0$ is very small and $b$ is very large here.

In 1961, Барбан [4] proved that equation (1.1) holds when $\eta_0 = 1/6$.

In 1962, Pan [5] proved that equation (1.1) holds when $\eta_0 = 1/3$, and obtained Proposition {1, 5}. Subsequently, Wang [6] obtained Proposition {1, 4} from $\eta_0 = 1/3$, and further deduced a not so apparent relation between $\eta_0$ and $b$: that Proposition {1, 4} and Proposition {1, 3} can be deduced from $\eta_0 = 1/3.327$ and $\eta_0 = 1/2.475$ respectively.

Pan [7] (1962) and Барбан [8] (1963) proved independently that equation (1.1) holds when $\eta_0 = 3/8$, and obtained Proposition {1, 4}. Бухштаб [9] obtained Proposition {1, 3} in 1965.

Almost at the same time, Виноградов [10] and Bombieri [11] independently proved that equation (1.1) holds when $\eta_0 = 1/2$, and also obtained Proposition {1, 3}.

In 1973, Chen [12] obtained Proposition {1, 2} by using sieves with weights.

The sieve methods used by the mathematicians above must sift out all the composite numbers. Although they have attained great success, up until now, Proposition {1, 1} is still unsolved by the sieve methods. In fact, it is not clear that the sieve methods that must sift out all the composite numbers will be able to prove Proposition {1, 1} ultimately. Indeed, if the ultimate aim is to prove Proposition {1, 1} rather than just sifting out all the composites in the Goldbach problem, it is immaterial whether all the composites are sifted out or not as long as the Proposition {1, 1} is proved.

Can Proposition {1, 1} be proved under the condition that only some composites are sifted out? The answer is positive. Indeed, one algebraic method for proving Proposition {1, 1} has already been found by Song [13]. This paper shows that three methods called *sifting function partition by integer sort*, *sifting function partition by intervals* and *comparative sieve method* respectively can also be used to prove Proposition {1, 1} without having to sift out all the composites. In fact, it is only necessary to show that there are prime numbers left in the residual integers after the initial sifting. These are feasible approachs for solving the Goldbach problem and can be practically implemented.

## 2. Some existing results of the Selberg's sieve

Let $\mathscr{A}$ be a nonempty set of positive integers not exceeding $N$, $\pi$ be the set of all primes, $\mathscr{P}$ be a subset of $\pi$, and $|\mu|$ be the number of elements of the set $\mu$. The aim of the Brun-Selberg sieve is to obtain bounds for the "sifting function" $S(\mathscr{A}; \mathscr{P}, z)$:

$$S(\mathscr{A};\mathscr{P},z) = \left|\{a | \ a \in \mathscr{A}, \ (a, P(z)) = 1\}\right|, \tag{2.1}$$

Where $z \geq 2$ is a real number and

$$P(z) = \prod_{\substack{p<z \\ p \in \mathscr{P}}} p. \tag{2.2}$$

Obviously, $S(\mathscr{A}; \mathscr{P}, z)$ gives the number of elements that are coprime with $P(z)$ in the set $\mathscr{A}$.

Let $\mathscr{A}_d$ be the set of those integers that can be divided by $d$ in the set $\mathscr{A}$. Then

$$S(\mathscr{A};\mathscr{P},z) = \sum_{d | P(z)} \mu(d) |\mathscr{A}_d|.$$

Clearly, $S(\mathscr{A}; \mathscr{P}, z)$ can be determined by calculating $|\mathscr{A}_d|$. In Selberg's sieve, one can choose a non-negative multiplicative function $\omega(d)$, and let $\dfrac{\omega(d)}{d}|\mathscr{A}|$ represent $|\mathscr{A}_d|$ approximatively, with an error term

$$r_d = |\mathscr{A}_d| - \frac{\omega(d)}{d}|\mathscr{A}|.$$

Then $S(\mathscr{A}; \mathscr{P}, z)$ can be determined by estimating the remainder term, i.e. the error term.

Accurate estimates of the upper and lower bounds for $S(\mathscr{A}; \mathscr{P}, z)$ were obtained by Jurkat-Richert [14] as follows (see also p.198 of [15]):

**Lemma 2.1.** *Suppose that $|\mathscr{A}| = X[1 + O(\ln^{-1} X)]$, where $X$ is a real number, $L_1 > 1$ is a constant. Suppose further that when $2 \leq z \leq X$, the multiplicative function $\omega(d)$ satisfies*

$$0 < \frac{\omega(p)}{p} \leq 1 - \frac{1}{L_1}, \quad p \in \mathscr{P}, \tag{2.3}$$

$$\sum_{w \leq p < z} \frac{\omega(p) \ln p}{p} = \ln \frac{z}{w} + O(1), \quad 2 \leq w < z; \tag{2.4}$$

*and that there exist real numbers $\alpha$ and $B$ satisfying $0 < \alpha \leq 1$ and $B > 0$ such that the following estimation of the remainder term holds*:

$$\sum_{\substack{d \mid P(z) \\ d \leq X^\alpha \ln^{-B} X}} \mu^2(d) 3^{v_1(d)} |r_d| \ll \frac{X}{\ln^2 X}. \tag{2.5}$$

*Then the upper and lower bounds for the sifting function $S(\mathscr{A}; \mathscr{P}, z)$ satisfy*

$$S(\mathscr{A};\mathscr{P},z) \leq XW(z) F\left(\alpha \frac{\ln X}{\ln z}\right)\left(1 + O\left(\frac{1}{\ln^{1/14} X}\right)\right), \tag{2.6}$$

$$S(\mathscr{A};\mathscr{P},z) \geq XW(z) f\left(\alpha \frac{\ln X}{\ln z}\right)\left(1 + O\left(\frac{1}{\ln^{1/14} X}\right)\right) \tag{2.7}$$

*respectively, where*

$$W(z) = \prod_{\substack{p < z \\ p \in \mathscr{P}}} \left(1 - \frac{\omega(p)}{p}\right), \tag{2.8}$$

*$F(u)$ is a decreasing function with $1 \leq F(u) \leq 2e^\gamma$ and $f(u)$ is an increasing function with $0 \leq f(u) \leq 1$ when $u \geq 1$, and*

$$F(u) = f(u) + O\left(\frac{e^{-u}}{u}\right) = 1 + O(e^{-u}). \tag{2.9}$$

The following formulas are due to Mertens (refer to [15]):

**Lemma 2.2.**
$$\prod_{p < z}\left(1 - \frac{1}{p}\right) = \frac{e^{-\gamma}}{\ln z} + O\left(\frac{1}{\ln^2 z}\right), \tag{2.10}$$

$$\sum_{p<z} \frac{\ln p}{p-1} = \sum_{p<z} \frac{\ln p}{p} + O(1) = \ln z + O(1). \tag{2.11}$$

**Lemma 2.3.** *Assume that $\mu(n)$ is the Möbius function, $\nu_1(n)$ is the number of different prime divisors of n and d(n) is the divisor function. Then*
$$\mu^2(n)3^{\nu_1(n)} \leq d^2(n).$$

*Proof.* When $\mu(n) = 0$, we have $\mu^2(n)3^{\nu_1(n)} = 0$, but $d(n) > 0$. Thus $\mu^2(n)3^{\nu_1(n)} < d^2(n)$, and the lemma holds.

When $n = 1$, we have $\mu(n) = 1$, $\nu_1(1) = 0$ and $d(1) = 1$. Thus $\mu^2(1)3^{\nu_1(1)} = d^2(1)$, and the lemma holds.

When $n > 1$ and $\mu(n) \neq 0$, then $\mu^2(n) = 1$ and $n$ is square-free. Thus $d(n) = 2^{\nu_1(n)}$ (which can be proved by induction), and $d^2(n) = 4^{\nu_1(n)}$. Thus $\mu^2(n)3^{\nu_1(n)} < d^2(n)$, and the lemma holds.

The following lemma appears as Corollary 1 of Theorem 1 in Chapter eight §1 of [15] (p.208):

**Lemma 2.4.** *Assume that $x \geq 2$ is a real number. The following estimate of the remainder term holds when $B = 38$:*
$$\sum_{d \leq x^{1/2} \ln^{-B} x} \mu^2(d)3^{\nu_1(d)} \max_{y \leq x} \max_{(l,d)=1} \left| \pi(y;l,d) - \frac{\operatorname{Li}(y)}{\varphi(d)} \right| \ll \frac{x}{\ln^3 x}.$$

The following lemma appears as Lemma 2 of Chapter three §1 in [15] (p.56):

**Lemma 2.5.** *Suppose that $x \geq 2$ is a real number, $r \geq 0$ is an integer, and $d(n)$ is the divisor function. Then*
$$\sum_{1 \leq n \leq x} d^r(n) \ll x(\ln x)^{2^r - 1}.$$

Suppose that $N$ is a large integer. Introduce a set $\mathscr{B}$ of integers and its subset $B(m)$:
$$\mathscr{B} = \{b \mid 1 < b \leq N\},$$
$$B(m) = \{b \mid b \in \mathscr{B};\ \text{if } p \mid b \text{ then } p \geq N^{1/m}\}.$$
Obviously, $|\mathscr{B}| = N - 1$ holds, and $|B(m)| = S(\mathscr{B}; \mathscr{P}, N^{1/m})$ can be determined by sieve methods.

**Lemma 2.6.** *Suppose that* $|r_d| = O(1)$, $D = N \ln^{-5} N$. *Then*

$$\sum_{d \leq D} \mu^2(d) 3^{\nu_1(d)} |r_d| \ll \frac{N}{\ln^2 N}.$$

*Proof.* Because of $|r_d| = O(1)$, from Lemma 2.3 one has $\mu^2(n) 3^{\nu_1(n)} \leq d^2(n)$. Then, from Lemma 2.5, one obtains

$$\sum_{d \leq D} \mu^2(d) 3^{\nu_1(d)} |r_d| \ll \sum_{n \leq D} d^2(n) \ll D (\ln D)^{2^2 - 1} \ll \frac{N}{\ln^2 N}.$$

The lemma holds.

**Theorem 2.7.** *Suppose that N is a large even integer and $m \geq 1$ is a real number. The following estimates hold*:

$$|B(m)| \leq m e^{-\gamma} \frac{N}{\ln N} F(m) \left(1 + O\left(\frac{1}{\ln^{1/14} N}\right)\right), \tag{2.12}$$

$$|B(m)| \geq m e^{-\gamma} \frac{N}{\ln N} f(m) \left(1 + O\left(\frac{1}{\ln^{1/14} N}\right)\right). \tag{2.13}$$

*Proof.* In order to sift out some composites in $\mathscr{B}$, put $\mathscr{P} = \pi$, $\omega(d) = 1$, $z = N^{1/m}$ and

$$\mathscr{B}_d = \{a | a \in \mathscr{B}, \ d | a, \ d | P(z)\}$$

in Lemma 2.1. Since $|\mathscr{B}| = N - 1$, let $X = N$; and then one obtains

$$r_d = |\mathscr{B}_d| - \frac{N}{d}, \ |r_d| < 1;$$

thus equation (2.5) holds when $\alpha = 1$ (see Lemma 2.6). Because the least prime in $\mathscr{P}$ is 2, $\omega(d)$ satisfies the conditions (2.3) and (2.4) when $L_1 = 2$. Therefore, from Lemma 2.2, we have

$$W(z) = \prod_{p < N^{1/m}} \left(1 - \frac{1}{p}\right) = \frac{m e^{-\gamma}}{\ln N} \left(1 + O\left(\frac{1}{\ln N}\right)\right).$$

Noting that $F(m) \geq 1$ and $f(m) \geq 0$ are well-defined when $m \geq 1$, one obtains, from Lemma 2.1 and noting that $S(\mathscr{B}; \mathscr{P}, N^{1/m}) = |B(m)|$,

$$|B(m)| \leq m e^{-\gamma} \frac{N}{\ln N} F(m) \left(1 + O\left(\frac{1}{\ln^{1/14} N}\right)\right),$$

$$|B(m)| \geq m e^{-\gamma} \frac{N}{\ln N} f(m) \left(1 + O\left(\frac{1}{\ln^{1/14} N}\right)\right).$$

The theorem is proved.

Suppose that $N$ is a large even integer. Define the set $\mathscr{A}$ and its subset $\mathscr{A} \cap \pi$ as follows:

$$\mathscr{A} = \{a \mid a = N - p,\ p \in \pi,\ p \leq N\},$$

$$\mathscr{A} \cap \pi = \{a \mid a \in \mathscr{A},\ a \in \pi\}.$$

Obviously, $\mathscr{A}$ is a subset of $\mathscr{B}$ and $|\mathscr{A}| = \pi(N)$. Furthermore, $|\mathscr{A} \cap \pi|$ is the number of representations of $N$ as the sum of two primes. This is precisely the answer to the Goldbach problem.

Next, we calculate the sifting function $S(\mathscr{A};\ \mathscr{P},\ N^{1/m}) = |\mathscr{A} \cap B(m)|$, where $\mathscr{A} \cap B(m)$ is a subset of the set $B(m)$.

**Theorem 2.8.** *Suppose that $N$ is a large even integer and $m \geq 2$ is a real number. The following estimates hold*:

$$|\mathscr{A} \cap B(m)| \leq 2me^{-\gamma} \frac{c(N)N}{\ln^2 N} F(m/2)\left(1 + O\left(\frac{1}{\ln^{1/14} N}\right)\right), \qquad (2.14)$$

$$|\mathscr{A} \cap B(m)| \geq 2me^{-\gamma} \frac{c(N)N}{\ln^2 N} f(m/2)\left(1 + O\left(\frac{1}{\ln^{1/14} N}\right)\right), \qquad (2.15)$$

*where*

$$c(N) = \prod_{\substack{p \mid N \\ p > 2}} \frac{p-1}{p-2} \prod_{p > 2}\left(1 - \frac{1}{(p-1)^2}\right). \qquad (2.16)$$

*Proof.* In order to sift out some composites in $\mathscr{A}$, put $\mathscr{P} = \{p \mid p \in \pi,\ (p, N) = 1\}$, $z = N^{1/m}$, $\omega(d) = \dfrac{d}{\varphi(d)}$, $(d, N) = 1$ and

$$\mathscr{A}_d = \{a \mid\ a \in \mathscr{A},\ d \mid a,\ d \mid P(z)\}$$

in Lemma 2.1. Since $|\mathscr{A}| = \pi(N)$, let $X = \mathrm{Li}(N)$; then one obtains

$$r_d = |\mathscr{A}_d| - \frac{\omega(d)}{d} X = \pi(N; N, d) - \frac{\mathrm{Li}(N)}{\varphi(d)}.$$

Thus, equation (2.5) holds when $\alpha = 0.5$ and $B = 38$ (see Lemma 2.4). Because $2 \mid N$, we see that the least prime in $\mathscr{P}$ is not less than 3 and so all the conditions of Lemma 2.1 are satisfied when $L_1 = 2$. Furthermore, we have

$$W(N^{1/m}) = \prod_{\substack{p\nmid N \\ p<N^{1/m}}} \left(1 - \frac{\omega(p)}{p}\right) = \prod_{\substack{p\nmid N \\ p<N^{1/m}}} \frac{p-2}{p-1}$$

$$= \prod_{\substack{p\mid N \\ p<N^{1/m}}} \frac{p}{p-1} \prod_{\substack{p\nmid N \\ p<N^{1/m}}} \left(1 - \frac{1}{(p-1)^2}\right) \prod_{p<N^{1/m}} \left(1 - \frac{1}{p}\right)$$

$$= \frac{2m\mathrm{e}^{-\gamma}}{\ln N} \prod_{\substack{p\mid N \\ 2<p<N^{1/m}}} \frac{p-1}{p-2} \prod_{2<p<N^{1/n}} \left(1 - \frac{1}{(p-1)^2}\right)\left(1 + O\left(\frac{1}{\ln N}\right)\right)$$

$$= \frac{2m\mathrm{e}^{-\gamma}}{\ln N} \prod_{\substack{p\mid N \\ p>2}} \frac{p-1}{p-2} \prod_{p>2} \left(1 - \frac{1}{(p-1)^2}\right)\left(1 + O\left(\frac{1}{\ln N}\right) + O\left(\frac{1}{N^{1/m}}\right)\right)$$

$$= 2m\mathrm{e}^{-\gamma} \frac{c(N)}{\ln N}\left(1 + O\left(\frac{1}{\ln N}\right)\right).$$

Noting that $F(m/2) \geq 1$ and $f(m/2) \geq 0$ are well-defined when $m \geq 2$, one obtains

$$S(\mathscr{A};\mathscr{P}, N^{1/m}) \leq 2m\mathrm{e}^{-\gamma} \frac{c(N)N}{\ln^2 N} F(m/2)\left(1 + O\left(\frac{1}{\ln^{1/14} N}\right)\right),$$

$$S(\mathscr{A};\mathscr{P}, N^{1/m}) \geq 2m\mathrm{e}^{-\gamma} \frac{c(N)N}{\ln^2 N} f(m/2)\left(1 + O\left(\frac{1}{\ln^{1/14} N}\right)\right),$$

where $S(\mathscr{A}; \mathscr{P}, N^{1/m}) = |\mathscr{A} \cap B(m)|$. Therefore, the theorem holds.

## 3. The algebraic relation between $|\mathscr{A}\cap B(m)|$ and $|B(m)|$

The sifting functions $|\mathscr{A} \cap B(m)|$ and $|B(m)|$ determined by sieve methods are shown in equations (2.12) ~ (2.15). These are clearly homogeneous linear correlation functions whose exact form can be determined by algebraic methods.

**Theorem 3.1.** *Suppose that N is a large even number and $m \geq 2$ is a real number. The following relation between the two sifting functions $|\mathscr{A}\cap B(m)|$ and $|B(m)|$ holds*:

$$|\mathscr{A} \cap B(m)| = 2\alpha(m,N)\frac{c(N)}{\ln N}|B(m)|\left(1 + O\left(\frac{1}{\ln^{1/14} N}\right)\right), \tag{3.1}$$

*where $\alpha(m,N)$ with $\frac{f(m/2)}{F(m)} \leq \alpha(m,N) \leq \frac{F(m/2)}{f(m)}$ and $\lim_{m\to\infty} \alpha(m,N) = 1$ is a real number.*

*Proof.* From equations (2.12) and (2.15) one obtains

$$\frac{|\mathscr{A} \cap B(m)|}{2\frac{c(N)}{\ln N}|B(m)|} \geq \frac{f(m/2)}{F(m)}\left(1+O\left(\frac{1}{\ln^{1/14} N}\right)\right).$$

From equations (2.13) and (2.14) one obtains

$$\frac{|\mathscr{A} \cap B(m)|}{2\frac{c(N)}{\ln N}|B(m)|} \leq \frac{F(m/2)}{f(m)}\left(1+O\left(\frac{1}{\ln^{1/14} N}\right)\right).$$

Therefore, there must exist a real number $\alpha(m, N)$ with $\frac{f(m/2)}{F(m)} \leq \alpha(m,N) \leq \frac{F(m/2)}{f(m)}$ such that

$$\frac{|\mathscr{A} \cap B(m)|}{2\frac{c(N)}{\ln N}|B(m)|} = \alpha(m,N)\left(1+O\left(\frac{1}{\ln^{1/14} N}\right)\right)$$

or

$$|\mathscr{A} \cap B(m)| = 2\alpha(m,N)\frac{c(N)}{\ln N}|B(m)|\left(1+O\left(\frac{1}{\ln^{1/14} N}\right)\right)$$

holds. On the other hand, from equation (2.9) one obtains

$$\lim_{m\to\infty} \alpha(m,N) \leq \lim_{m\to\infty} \frac{F(m/2)}{f(m)} = 1,$$

$$\lim_{m\to\infty} \alpha(m,N) \geq \lim_{m\to\infty} \frac{f(m/2)}{F(m)} = 1.$$

Therefore the theorem holds.

The proof of Theorem 3.1 requires only algebraic operations and requires none of the techniques of the sieve methods. Intuitively, the value of $\alpha(m, N)$ appears to be related to the parameter $m$ because both the upper and lower bounds for $\alpha(m, N)$ are dependent on $m$. However, there could be exceptions; namely $\alpha(m, N)$ could possibly be a constant between both bounds. If it is so, then we must have $\alpha(m, N) = 1$ because this agrees with $\lim_{m\to\infty} \alpha(m, N) = 1$. In fact, it can be proved that $\alpha(m, N)$ is exactly so, see below.

The homogeneous linear relation between $|\mathscr{A} \cap B(m)|$ and $|B(m)|$ given by Theorem 3.1 is related to the fact that $\mathscr{A} \cap B(m)$ is a subset of $B(m)$. For obviously if $|B(m)| = 0$, then we must have $|\mathscr{A} \cap B(m)| = 0$ which can only be true if the relation between $|\mathscr{A} \cap B(m)|$ and $|B(m)|$ is homogeneous. Furthermore, if $\mathscr{A} \cap B(m)$ is a subset of $B(m)$ and the integers within $\mathscr{A} \cap B(m)$ are randomly distributed, (in other words when the distribution of primes in the set of integers is random), then the relation between $|\mathscr{A} \cap B(m)|$ and $|B(m)|$ should be linear. In effect, Theorem 3.1 shows that the distribution of primes in the set of integers is random.

## 4. Sifting function partition by integer sort and Proposition {1, 1}

Although Proposition {1, 1} cannot be proved by using any of the existing sieve methods alone, the sifting functions given by these methods can be used in combination with algebraic techniques to solve both the Goldbach problem and the problem of twin primes. Sifting functions generate estimates of the number of integers in certain integer sets, and integer sets contain information about the various properties of integers. Interestingly, this information can only be revealed by algebraic methods and not by sieve methods because only algebraic methods are able to separate those integers with specific properties from the set. By determining the cardinality of these integers, using sifting functions, some problems that cannot be solved solely by means of the sieve methods can possibly be solved. See below.

In order to prove Proposition {1, 1}, classify the integers according to the number of their prime divisors. Suppose that $k \geq 1$ is an integer, let $B_k(m)$ be a set of integers with $k$ prime divisors:

$$B_k(m) = \{b_k \mid b_k = p_1 \cdots p_k, \ N^{1/m} \leq p_k \leq \cdots \leq p_1, \ b_k \leq N\}.$$

Obviously, $B_k(m)$ is a subset of $B(m)$, and the following relations hold:

$$\mathscr{A} \cap B_k(m) = B_k(m) = \varnothing, \ \left|\mathscr{A} \cap B_k(m)\right| = \left|B_k(m)\right| = 0 \ \text{if} \ k > m,$$

$$B_i(m) \cap B_j(m) = \varnothing \ \text{if} \ i \neq j,$$

$$B(m) = \bigcup_{k \geq 1} B_k(m) \ \text{and}$$

$$\mathscr{A} \cap B(m) = \bigcup_{k \geq 1} (\mathscr{A} \cap B_k(m)).$$

The following partitions of the sifting functions $|\mathscr{A} \cap B(m)|$ and $|B(m)|$ hold:

$$\left|\mathscr{A} \cap B(m)\right| = \sum_{k \geq 1} \left|\mathscr{A} \cap B_k(m)\right|,$$

$$\left|B(m)\right| = \sum_{k \geq 1} \left|B_k(m)\right|.$$

Now we deduce the relation between the *partitions of the sifting functions* $|\mathscr{A} \cap B(m)|$ *and* $|B(m)|$.

**Theorem 4.1.** *Suppose that N is a large even number, $m \geq 2$ is a real number and $k \geq 1$ is an integer. The following relation between the partitions of the sifting functions $|\mathscr{A} \cap B(m)|$ and $|B(m)|$ holds*:

$$|\mathscr{A} \cap B_k(m)| = 2\alpha(m,N)\frac{c(N)}{\ln N}|B_k(m)|\left(1+O\left(\frac{1}{\ln^{1/14} N}\right)\right). \qquad (4.1)$$

*Proof.* The theorem can be proved by means of the "separation of variables". Introduce the subset $C_k(m)$ of the set $B(m)$:

$$C_k(m) = \bigcup_{\substack{i \geq 1 \\ i \neq k}} B_i(m).$$

Obviously, the following relations hold:

$$B_k(m) \cap C_k(m) = \emptyset,$$

$$B(m) = B_k(m) \cup C_k(m),$$

$$\mathscr{A} \cap B(m) = (\mathscr{A} \cap B_k(m)) \cup (\mathscr{A} \cap C_k(m)).$$

Also, the following partitions of the sifting functions $|\mathscr{A} \cap B(m)|$ and $|B(m)|$ hold:

$$|B(m)| = |B_k(m)| + |C_k(m)|, \qquad (4.2)$$

$$|\mathscr{A} \cap B(m)| = |\mathscr{A} \cap B_k(m)| + |\mathscr{A} \cap C_k(m)|. \qquad (4.3)$$

From Theorem 3.1 and equations (4.2) and (4.3) one obtains

$$|\mathscr{A} \cap B_k(m)| + |\mathscr{A} \cap C_k(m)| =$$

$$= 2\alpha(m,N)\frac{c(N)}{\ln N}\left(1+O\left(\frac{1}{\ln^{1/14} N}\right)\right)(|B_k(m)| + |C_k(m)|), \qquad (4.4)$$

where $|\mathscr{A} \cap B_k(m)|$ and $|\mathscr{A} \cap C_k(m)|$ are undetermined functions, which can be determined by using the method of "*separation of variables*". In fact, equation (4.4) can be written in the form

$$|\mathscr{A} \cap B_k(m)| - 2\alpha(m,N)\frac{c(N)}{\ln N}|B_k(m)|\left(1+O\left(\frac{1}{\ln^{1/14} N}\right)\right) =$$

$$= 2\alpha(m,N)\frac{c(N)}{\ln N}|C_k(m)|\left(1+O\left(\frac{1}{\ln^{1/14} N}\right)\right) - |\mathscr{A} \cap C_k(m)|. \qquad (4.5)$$

Where the constants of the two big O notations may be different. Let

$$c_1 = \frac{2\alpha(m,N)c(N)}{\ln N}\left(1+O\left(\frac{1}{\ln^{1/14} N}\right)\right),$$

$$c_2 = \frac{2\alpha(m,N)c(N)}{\ln N}\left(1+O\left(\frac{1}{\ln^{1/14} N}\right)\right),$$

one obtains

$$|\mathscr{A} \cap B_k(m)| - c_1|B_k(m)| = c_2|C_k(m)| - |\mathscr{A} \cap C_k(m)| = F,$$

where $F$ is an undetermined function.

Note that

$$C_k(m) \cap B_k(m) = \varnothing,$$
$$C_k(m) \cap [\mathscr{A} \cap B_k(m)] = \varnothing,$$
$$[\mathscr{A} \cap C_k(m)] \cap B_k(m) = \varnothing,$$
$$[\mathscr{A} \cap C_k(m)] \cap [\mathscr{A} \cap B_k(m)] = \varnothing.$$

From $F = |\mathscr{A} \cap B_k(m)| - c_1|B_k(m)|$ we deduce that $F$ is only possibly related to $|B_k(m)|$ and $|\mathscr{A} \cap B_k(m)|$, but it is definitely independent from $|C_k(m)|$ and $|\mathscr{A} \cap C_k(m)|$. Likewise, from $F = c_2|C_k(m)| - |\mathscr{A} \cap C_k(m)|$ we deduce that $F$ is also definitely independent from $|B_k(m)|$ and $|\mathscr{A} \cap B_k(m)|$.

To sum up, therefore, $F$ is definitely independent from $|B_k(m)|$, $|\mathscr{A} \cap B_k(m)|$, $|C_k(m)|$ and $|\mathscr{A} \cap C_k(m)|$, i.e. $F$ is a constant independent from the parameters $k$ and $m$, although it can still possibly be related to $|\mathscr{A}|$, i.e. $F = F(|\mathscr{A}|)$. And then from equation (4.5) one obtains

$$|\mathscr{A} \cap B_k(m)| - 2\alpha(m,N)\frac{c(N)}{\ln N}|B_k(m)|\left(1+O\left(\frac{1}{\ln^{1/14} N}\right)\right) = F(|\mathscr{A}|),$$

$$|\mathscr{A} \cap C_k(m)| - 2\alpha(m,N)\frac{c(N)}{\ln N}|C_k(m)|\left(1+O\left(\frac{1}{\ln^{1/14} N}\right)\right) = -F(|\mathscr{A}|).$$

Since $F = F(|\mathscr{A}|)$ is independent from the parameters $k$ and $m$, noting that, in the particular case that when $k > m$, we have $|B_k(m)| = 0$, $|\mathscr{A} \cap B_k(m)| = 0$, $|C_k(m)| = |B(m)|$ and $|\mathscr{A} \cap C_k(m)| = |\mathscr{A} \cap B(m)|$, simultaneously; and then $F(|\mathscr{A}|) = 0$ holds for any $|B_k(m)|$, $|\mathscr{A} \cap B_k(m)|$, $|C_k(m)|$ and $|\mathscr{A} \cap C_k(m)|$ with any $k$ and $m$. The theorem holds.

The proof of Proposition {1, 1} can now be completed and the indeterminacy of $\alpha(m, N)$ can also be eliminated in one stroke because it can be proved that we must have $\alpha(m, N) = 1$:

**Theorem 4.2.** *Suppose that N is a large even number. The following estimate holds*:

$$|\mathscr{A} \cap \pi| = 2c(N)\frac{N}{\ln^2 N}\left(1 + O\left(\frac{1}{\ln^{1/14} N}\right)\right). \tag{4.6}$$

*Proof.* When $k = 1$, one obtains, from Theorem 4.1, that

$$|\mathscr{A} \cap B_1(m)| = 2\alpha(m, N)\frac{c(N)}{\ln N}|B_1(m)|\left(1 + O\left(\frac{1}{\ln^{1/14} N}\right)\right).$$

Since

$$|\mathscr{A} \cap \pi| = |\mathscr{A} \cap B_1(m)| + O(\pi(N^{1/m})),$$

$$|B_1(m)| = \pi(N) - \pi(N^{1/m}) = \frac{N}{\ln N}\left(1 + O\left(\frac{1}{\ln N}\right)\right),$$

we must have

$$|\mathscr{A} \cap \pi| = 2\alpha(m, N)\frac{c(N)N}{\ln^2 N}\left(1 + O\left(\frac{1}{\ln^{1/14} N}\right)\right).$$

However, because $|\mathscr{A} \cap \pi|$ is the number of representations of $N$ as the sum of two primes, whose value is independent from the parameter $m$, therefore, $\alpha(m, N)$ must be independent from $m$. Noting that $\lim_{m \to \infty} \alpha(m, N) = 1$, we must have $\alpha(m, N) = 1$ for any real number $m \geq 2$. Thus the theorem holds.

Theorem 4.2 is the final answer to the Goldbach problem! This problem, which cannot be solved by using the sieve methods alone, can now be satisfactorily solved by using a combination of both sieve methods and algebraic methods in the manner discussed above.

An even more accurate and general result to the Goldbach problem was obtained by Song in 1999 [16]. By using probability methods, it was proved that the number (expectation) of representations of any positive integer $N$ as a sum of $n$ odd primes is equal to

$$D_n(N) = \sum_{\substack{1<x_i<N \\ 1\leq i\leq n-1}} \left( C_n(x_n) \prod_{1\leq i\leq n} r(x_i) \right), \qquad (4.7)$$

where $x_n = N - (x_1 + \ldots + x_{n-1})$;

$$C_n(x_n) = \prod_{\substack{p|N \\ p\leq x_n^c}} \left(1 + \frac{(-1)^n}{(p-1)^{n-1}}\right) \prod_{\substack{p\nmid N \\ p\leq x_n^c}} \left(1 - \frac{(-1)^n}{(p-1)^n}\right),$$

$c = e^{-\gamma}$, $\gamma$ is the Euler constant; and

$$r(x) = \begin{cases} 0, & x < 2, \\ \dfrac{1}{x} \sum_{k\geq 1} \dfrac{1}{\zeta(k+1)} \dfrac{\ln^{k-1} x}{k!}, & x \geq 2, \end{cases}$$

where $\zeta(k)$ is the Riemann's zeta function. Note that when $2 \nmid n+N$, we have $C_n(x_n) = 0$ and $D_n(N) = 0$, implying that $D_n(N)$ can only be non-zero when the parities of $n$ and $N$ are the same.

The asymptotic formula of equation (4.7) is [17]

$$D_n(N) \sim \frac{c_n(N)}{(n-1)!} \frac{N^{n-1}}{\ln^n N}, \qquad (4.8)$$

where

$$c_n(N) = \prod_{p|N} \left(1 + \frac{(-1)^n}{(p-1)^{n-1}}\right) \prod_{p\nmid N} \left(1 - \frac{(-1)^n}{(p-1)^n}\right).$$

It can be easily proved that $c_2(N) = 2c(N)$ and equation (4.8) is the same as equation (4.6) when $n = 2$ and $2|N$. However equation (4.8) is applicable to any even integers with $N \geq 6$, and equation (4.7) yields far more precise results than equation (4.6). When $n = 3$ and $N \geq 9$ ($2 \nmid N$), we have

$$c_3(N) = \prod_{p|N} \left(1 - \frac{1}{(p-1)^2}\right) \prod_{p\nmid N} \left(1 + \frac{1}{(p-1)^3}\right),$$

so that equation (4.8) is consistent with Виноградов's result [18] in this case. Similarly, equation (4.8) is also applicable to small odd numbers $N$, and equation (4.7) yields far more precise results than Виноградов's result.

## 5. Sifting function partition by intervals and the Proposition {1, 1}

Although Proposition {1, 1} cannot be proved by using any of the existing sieve methods alone, the sifting functions given by these methods can be used in combination with algebraic techniques to solve both the Goldbach problem and the problem of twin primes. Sifting functions generate estimates of the number of integers in certain integer sets, and integer sets contain information about the various properties of integers. Interestingly, this information can only be revealed by algebraic methods and not by sieve methods because only algebraic methods are able to separate those integers with specific properties from the set. By determining the cardinality of these integers, using sifting functions, some problems that cannot be solved solely by means of the sieve methods can possibly be solved. See below.

A simple method for proving Proposition $\{1, 1\}$ is to separate those primes not exceeding $N^{2/m}$ from the sets $\mathcal{A} \cap B(m)$ and $B(m)$. Firstly, partition the set $B(m)$ by intervals into two subsets $B_l(m, w)$ and $B_r(m, w)$:

$$B_l(m, w) = \{b \mid b \in B(m),\ N^{1/m} \leq b \leq w\},$$

$$B_r(m, w) = \{b \mid b \in B(m),\ b > w\},$$

where $w \leq 0.5\,N$ is a real number. The following relations hold:

$$B_l(m, w) = \varnothing \text{ and } |B_l(m, w)| = 0 \text{ if } w \leq N^{1/m};$$

$$B_l(m, w) \cap B_r(m, w) = \varnothing;$$

$$B(m) = B_l(m, w) \cup B_r(m, w).$$

Obviously, the following *partition by intervals* of the sifting function $|B(m)|$ also holds:

$$|B(m)| = |B_l(m, w)| + |B_r(m, w)|. \tag{5.1}$$

Secondly, partition the set $\mathcal{A} \cap B(m)$ by intervals into two subsets $\mathcal{A} \cap B_l(m, w)$ and $\mathcal{A} \cap B_r(m, w)$:

$$\mathcal{A} \cap B(m) = [\mathcal{A} \cap B_l(m, w)] \cup [\mathcal{A} \cap B_r(m, w)].$$

The following relations hold:

$$\mathcal{A} \cap B_l(m, w) = \varnothing \text{ and } |\mathcal{A} \cap B_l(m, w)| = 0 \text{ if } w \leq N^{1/m};$$

$$[\mathcal{A} \cap B_l(m, w)] \cap [\mathcal{A} \cap B_r(m, w)] = \varnothing.$$

Obviously, the following *partition by intervals* of the sifting function $|\mathcal{A} \cap B(m)|$ also holds:

$$|\mathscr{A} \cap B(m)| = |\mathscr{A} \cap B_l(m,w)| + |\mathscr{A} \cap B_r(m,w)|. \tag{5.2}$$

Now deduce the relation of the partitions of the sifting functions $|\mathscr{A} \cap B(m)|$ and $|B(m)|$.

**Theorem 5.1.** *Suppose that $N$ is a large even number, $m \geq 2$ is a real number. The following relation between the partitions of the sifting functions $|\mathscr{A} \cap B(m)|$ and $|B(m)|$ holds*:

$$\left. \begin{array}{l} |\mathscr{A} \cap B_l(m,w)| = 2\alpha(m,N)\dfrac{c(N)}{\ln N}|B_l(m,w)|\left(1+O\left(\dfrac{1}{\ln^{1/14} N}\right)\right), \\[2mm] |\mathscr{A} \cap B_r(m,w)| = 2\alpha(m,N)\dfrac{c(N)}{\ln N}|B_r(m,w)|\left(1+O\left(\dfrac{1}{\ln^{1/14} N}\right)\right). \end{array} \right\} \tag{5.3}$$

*Proof.* From equations (5.1) and (5.2), equation (3.1) can be written in

$$|\mathscr{A} \cap B_l(m,w)| + |\mathscr{A} \cap B_r(m,w)| =$$

$$= 2\alpha(m,N)\frac{c(N)}{\ln N}(|B_l(m,w)|+|B_r(m,w)|)\left(1+O\left(\frac{1}{\ln^{1/14} N}\right)\right), \tag{5.4}$$

where $|\mathscr{A} \cap B_l(m)|$ and $|\mathscr{A} \cap B_r(m)|$ are undetermined functions. Equation (5.4) can be written as

$$|\mathscr{A} \cap B_l(m,w)| - 2\alpha(m,N)\frac{c(N)}{\ln N}|B_l(m,w)|\left(1+O\left(\frac{1}{\ln^{1/14} N}\right)\right) =$$

$$= 2\alpha(m,N)\frac{c(N)}{\ln N}|B_r(m,w)|\left(1+O\left(\frac{1}{\ln^{1/14} N}\right)\right) - |\mathscr{A} \cap B_r(m,w)|. \tag{5.5}$$

Where, the constants of two big O notations in above are possibly different. Let

$$c_1 = 2\alpha(m,N)\frac{c(N)}{\ln N}\left(1+O\left(\frac{1}{\ln^{1/14} N}\right)\right),$$

$$c_2 = 2\alpha(m,N)\frac{c(N)}{\ln N}\left(1+O\left(\frac{1}{\ln^{1/14} N}\right)\right).$$

From equation (5.5) one obtains

$$|\mathscr{A} \cap B_l(m,w)| - c_1|B_l(m,w)| = c_2|B_r(m,w)| - |\mathscr{A} \cap B_r(m,w)| = F,$$

where $F$ is an undetermined function.

Note that

$$B_l(m, w) \cap B_r(m, w) = \emptyset,$$
$$B_l(m, w) \cap [\mathscr{A} \cap B_r(m, w)] = \emptyset,$$
$$[\mathscr{A} \cap B_l(m, w)] \cap B_r(m, w) = \emptyset,$$
$$[\mathscr{A} \cap B_l(m, w)] \cap [\mathscr{A} \cap B_r(m, w)] = \emptyset.$$

From $F = |\mathscr{A} \cap B_l(m, w)| - c_1 |B_l(m, w)|$ we deduce that $F$ is only possibly related to $|B_l(m,w)|$ and $|\mathscr{A} \cap B_l(m, w)|$, but it is definitely independent from $|B_r(m, w)|$ and $|\mathscr{A} \cap B_r(m, w)|$. Likewise, from $F = c_2 |B_r(m, w)| - |\mathscr{A} \cap B_r(m, w)|$ we deduce that $F$ is also definitely independent from $|B_l(m, w)|$ and $|\mathscr{A} \cap B_l(m, w)|$.

To sum up, therefore, $F$ is definitely independent from $|B_l(m,w)|$, $|\mathscr{A} \cap B_l(m, w)|$, $|B_r(m, w)|$ and $|\mathscr{A} \cap B_r(m, w)|$, i.e. $F$ is a constant independent from the parameters $m$ and $w$, although it can still possibly be related to $|\mathscr{A}|$, i.e. $F = F(|\mathscr{A}|)$. And then from equation (5.5) one obtains

$$|\mathscr{A} \cap B_l(m, w)| - 2\alpha(m, N) \frac{c(N)}{\ln N} |B_l(m, w)| \left(1 + O\left(\frac{1}{\ln^{1/14} N}\right)\right) = F(|\mathscr{A}|),$$

$$|\mathscr{A} \cap B_r(m, w)| - 2\alpha(m, N) \frac{c(N)}{\ln N} |B_r(m, w)| \left(1 + O\left(\frac{1}{\ln^{1/14} N}\right)\right) = -F(|\mathscr{A}|).$$

Since $F = F(|\mathscr{A}|)$ is independent from the parameters $m$ and $w$, noting that, in the particular case that $w \leq N^{1/m}$, we have $|B_l(m, w)| = 0$, $|\mathscr{A} \cap B_l(m, w)| = 0$, $|B_r(m, w)| = |B(m)|$ and $|\mathscr{A} \cap B_r(m, w)| = |\mathscr{A} \cap B(m)|$, simultaneously; and then $F(|\mathscr{A}|) = 0$ holds for any $|B_l(m, w)|$, $|B_l(m, w)|$, $|\mathscr{A} \cap B_r(m, w)|$ and $|B_r(m, w)|$ with any $m$ and $w$. The theorem holds.

The proof of the Proposition {1, 1} is complete.

**Theorem 5.2.** *The Proposition {1, 1} would be true if any Proposition {1, b} is true, and there is*

$$D(m, N) = m\alpha(m, N) c(N) \frac{N^{2/m}}{\ln^2 N} \left(1 + O\left(\frac{1}{\ln^{1/14} N}\right)\right). \qquad (5.6)$$

*Where $D(m, N)$ is the number of the primes not exceeding $N^{2/m}$ in the set $\mathscr{A}$, $m \geq 2$ and $\alpha(m, N)$ are real numbers with $\frac{f(m/2)}{F(m)} \leq \alpha(m, N) \leq \frac{F(m/2)}{f(m)}$, $b = [m]$ (or $b = m - 1$ if $m$ is an integer) is an integer.*

*Proof.* Put $w = N^{2/m}$ in the sets $B_l(m, w)$ and $B_r(m, w)$ of equation (5.3), noting that

$$|B_l(m, N^{2/m})| = \pi(N^{2/m}) - \pi(N^{1/m})$$

$$= \frac{mN^{2/m}}{2\ln N}\left(1 + O\left(\frac{1}{\ln N}\right)\right)$$

because the integers in the set $B_l(m, w)$ are all prime when $w = N^{2/m}$; from Theorem 5.1 one obtains

$$D(m, N) = |\mathscr{A} \cap B_l(m, N^{2/m})| + O(\pi(N^{1/m}))$$

$$= 2\alpha(m, N)\frac{c(N)}{\ln N}|B_l(m, N^{2/m})|\left(1 + O\left(\frac{1}{\ln^{1/14} N}\right)\right)$$

$$= m\alpha(m, N)c(N)\frac{N^{2/m}}{\ln^2 N}\left(1 + O\left(\frac{1}{\ln^{1/14} N}\right)\right).$$

The theorem holds.

The significance of Theorem 5.1 lies in the fact that as long as Proposition $\{1, b\}$ can be proved to be true for any integer $b$, then Proposition $\{1, 1\}$ will follow. For example, $f(2.5) = 0.8\, e^\gamma \ln 1.5$ when $m = 5$. From equation (2.15) one obtains

$$|\mathscr{A} \cap B(5)| \geq 8\ln 1.5 \frac{c(N)N}{\ln^2 N}\left(1 + O\left(\frac{1}{\ln^{1/14} N}\right)\right) > 0.$$

Therefore, Proposition $\{1, 4\}$ is true. Furthermore, since $\alpha(5, N) \geq \frac{f(5/2)}{F(5)} = 0.5767$ when $m = 5$, it follows from Theorem 5.2 that

$$D(5, N) \geq 2.8835\frac{c(N)N^{2/5}}{\ln^2 N}\left(1 + O\left(\frac{1}{\ln^{1/14} N}\right)\right) > 0.$$

Thus Proposition $\{1, 1\}$ also holds.

The uncertainty in the above results produced by the unknown value of the coefficient $\alpha(m, N)$ can be eliminated by proving that $\alpha(m, N) = 1$ for any $m \geq 2$, see below.

**Theorem 5.3.** *Suppose that N is a large even number. The following equation holds for any real number $m \geq 2$:*

$$\alpha(m, N) = 1. \tag{5.7}$$

*Proof.* When $m \geq 2$, Theorem 5.2 gives

$$D(m, N) = m\alpha(m, N)c(N)\frac{N^{2/m}}{\ln^2 N}\left(1 + O\left(\frac{1}{\ln^{1/14} N}\right)\right). \tag{5.8}$$

On the other hand, when $m' \geq 2$ with $\frac{m}{2} < m' < m$, Theorem 5.1 gives

$$|\mathscr{A} \cap B_l(m', w)| = 2\alpha(m', N)\frac{c(N)}{\ln N}|B_l(m', w)|\left(1 + O\left(\frac{1}{\ln^{1/14} N}\right)\right). \tag{5.9}$$

Put $w = N^{2/m}$ in the sets $\mathscr{A} \cap B_l(m', w)$ and $B_l(m', w)$ of equation (5.9), and calculating $D(m, N)$ from above expression. Noting that $\frac{2}{m'} > \frac{2}{m} > \frac{1}{m'}$ holds when $m' \geq 2$ with $\frac{m}{2} < m' < m$, which means that the integers within the sets $\mathscr{A} \cap B_l(m', w)$ and $B_l(m', w)$ and not exceeding $N^{2/m}$ (noting that $N^{2/m} < N^{2/m'}$) are all primes, thereout one obtains

$$|B_l(m', N^{2/m})| = \pi(N^{2/m}) - \pi(N^{1/m'})$$

$$= \frac{mN^{2/m}}{2\ln N}\left(1 + O\left(\frac{1}{\ln N}\right)\right),$$

Then, from Theorem 5.2, one obtains

$$D(m, N) = |\mathscr{A} \cap B_l(m', N^{2/m})| + O(\pi(N^{1/m'}))$$

$$= 2\alpha(m', N)\frac{c(N)}{\ln N}|B_l(m', N^{2/m})|\left(1 + O\left(\frac{1}{\ln^{1/14} N}\right)\right)$$

$$= m\alpha(m', N)\frac{c(N)N^{2/m}}{\ln^2 N}\left(1 + O\left(\frac{1}{\ln^{1/14} N}\right)\right). \tag{5.10}$$

It is clear, by comparing the $D(m, N)$ in equation (5.10) with that in equation (5.8), that $\alpha(m', N) = \alpha(m, N)$ since $D(m, N)$ and its first approximation are the same in both cases. The fact that $\alpha(m', N) = \alpha(m, N)$ must hold for any even integer $N$ and any real number pair $m \geq 2$ and $m' \geq 2$ with $\frac{m}{2} < m' < m$ implies that the value of $\alpha(m, N)$ is definitely not related to $m$. It is well known that $\lim_{m \to \infty} \alpha(m, N) = 1$, and so we must have $\alpha(m, N) = 1$ for any real number $m$. The theorem holds.

A first approximation to the Goldbach problem can be obtained as follows:

**Theorem 5.4.** *Suppose that $N$ is a large even number. The following estimate holds*:

$$|\mathscr{A} \cap \pi| = 2c(N)\frac{N}{\ln^2 N}\left(1 + O\left(\frac{1}{\ln^{1/14} N}\right)\right), \tag{5.11}$$

where $|\mathscr{A} \cap \pi|$ is the number of representations of N as the sum of two primes.

*Proof.* When $m = 2$ and $\alpha(m, N) = 1$, Theorem 3.1 gives

$$|\mathscr{A} \cap B(2)| = 2\frac{c(N)}{\ln N}|B(2)|\left(1 + O\left(\frac{1}{\ln^{1/14} N}\right)\right).$$

Noting that

$$|\mathscr{A} \cap \pi| = |\mathscr{A} \cap B(2)| + O(\pi(N^{1/2})),$$

$$|B(2)| = \pi(N) - \pi(N^{1/2}) = \frac{N}{\ln N}\left(1 + O\left(\frac{1}{\ln N}\right)\right),$$

Clearly, the theorem holds.

Theorem 5.4 is the best answer to the Goldbach problem! This problem which cannot be solved by using the sieve methods alone, can now be satisfactorily solved by using a combination of both sieve methods and algebraic methods in the manner discussed above.

## 6. The comparative sieve method and the Proposition {1, 1}

Suppose that $w < 0.5 N$ is a real number. In order to use the comparative sieve method, introduce a compared subset $\mathscr{B}_r(w)$ of $\mathscr{B}$ and its subset $B_r(m, w)$:

$$\mathscr{B}_r(w) = \{b \,|\, w < b \le N\},$$

$$B_r(m, w) = \{b \,|\, b \in \mathscr{B}_r(w);\ \text{if}\ p\,|\,b\ \text{then}\ p \ge N^{1/m}\}.$$

The sifting function $|B_r(m, w)| = S(\mathscr{B}_r(w);\ \mathscr{P},\ N^{1/m})$ can also be determined by using the sieve methods. Indeed, since $|\mathscr{B}_r(w)| = N - w$, let $X = N - w$, noting that

$$\sum_{d \le D} \mu^2(d) 3^{\nu_1(d)} |r_d| \ll \sum_{n \le D} d^2(n) \ll D(\ln D)^{2^2 - 1} \ll \frac{N}{\ln^2 N} \ll \frac{N - w}{\ln^2 N}$$

holds when $|r_d| = O(1)$, $D = N \ln^{-5} N$ and $w < 0.5 N$; as the proof of Theorem 2.7, the following theorem can be proved:

**Theorem 6.1.** *Suppose that N is a large integer, $w < 0.5 N$ and $m \ge 1$ are real numbers. The following estimates hold*:

$$|B_r(m, w)| \le me^{-\gamma}\frac{N - w}{\ln N}F(m)\left(1 + O\left(\frac{1}{\ln^{1/14} N}\right)\right), \tag{6.1}$$

$$|B_r(m,w)| \geq m e^{-\gamma} \frac{N-w}{\ln N} f(m)\left(1+O\left(\frac{1}{\ln^{1/14} N}\right)\right). \quad (6.2)$$

Moreover, define the compared subset $\mathscr{A}_r(w)$ of $\mathscr{A}$:

$$\mathscr{A}_r(w) = \{a \mid a = N-p,\ p \in \pi,\ p \leq N-w\},$$

where $w < 0.5N$ is a real number. Noting that $|\mathscr{A} \cap B_r(m, w)| = S(\mathscr{A}_r(w); \mathscr{P}, N^{1/m})$ can also be determined by using the sieve methods. Indeed, due to

$$|\mathscr{A}_r(w)| = \pi(N-w) = \frac{N-w}{\ln(N-w)}\left(1+O\left(\frac{1}{\ln(N-w)}\right)\right)$$

$$= \frac{N-w}{\ln N}\left(1+O\left(\frac{1}{\ln N}\right)\right),$$

Let $X = \mathrm{Li}(N-w)$, as the proof of Theorem 2.8, the following theorem can be proved:

**Theorem 6.2.** *Suppose that $N$ is a large even number, $w < 0.5N$ and $m \geq 2$ are real numbers. The following estimates hold*:

$$|\mathscr{A} \cap B_r(m,w)| \leq 2m e^{-\gamma} c(N) \frac{N-w}{\ln^2 N} F(m/2)\left(1+O\left(\frac{1}{\ln^{1/14} N}\right)\right), \quad (6.3)$$

$$|\mathscr{A} \cap B_r(m,w)| \geq 2m e^{-\gamma} c(N) \frac{N-w}{\ln^2 N} f(m/2)\left(1+O\left(\frac{1}{\ln^{1/14} N}\right)\right), \quad (6.4)$$

*Where $c(N)$ see also equation* (2.16).

Likewise, the relation between both compared sifting functions $|\mathscr{A} \cap B_r(m, w)|$ and $|B_r(m, w)|$ can also be determined from equations (6.1), (6.2), (6.3) and (6.4):

**Theorem 6.3.** *Suppose that $N$ is a large even number, $w < 0.5N$ and $m \geq 2$ are real numbers. The following relation between both sifting functions $|\mathscr{A} \cap B_r(m, w)|$ and $|B_r(m, w)|$ holds*:

$$|\mathscr{A} \cap B_r(m,w)| = 2\alpha_1(m,N)\frac{c(N)}{\ln N}|B_r(m,w)|\left(1+O\left(\frac{1}{\ln^{1/14} N}\right)\right), \quad (6.5)$$

here $\alpha_1(m, N)$ with $\frac{f(m/2)}{F(m)} \leq \alpha_1(m, N) \leq \frac{F(m/2)}{f(m)}$ and $\lim_{m \to \infty} \alpha_1(m, N) = 1$ is a real number.

Proposition {1, 1} cannot be proved by using the sieve methods because those composites whose divisors are all large primes can hardly be sifted out. However, the difficulty can be overcome in the comparative sieve method.

The key of the technology in the comparative sieve method is in that the set $\mathscr{A}$ should be partitioned into two subsets $\mathscr{A}_r(w)$ and $\mathscr{A}_l(w)$ with $\mathscr{A}_r(w) \cap \mathscr{A}_l(w) = \varnothing$ before sifting, and so do the sets $\mathscr{B}$, $\mathscr{B}_r(w)$ and $\mathscr{B}_l(w)$. Then we sift $\mathscr{A}$, $\mathscr{A}_r(w)$, $\mathscr{B}$ and $\mathscr{B}_r(w)$ by using the sieve method, which are completed in above. And compare both relations (3.1) and (3.2) at the last step.

Firstly, partition the set $B(m)$ into two subsets $B_l(m, w)$ and $B_r(m, w)$:
$$B_l(m, w) = \{b \mid b \in B(m),\ b \le w\},$$
$$B_r(m, w) = \{b \mid b \in B(m),\ b > w\},$$
where $w < 0.5N$ is a real number. Due to $B_l(m, w) \cap B_r(m, w) = \varnothing$, $B(m) = B_l(m, w) \cup B_r(m, w)$, the following partition of the sifting function $|B(m)|$ holds:
$$|B(m)| = |B_l(m, w)| + |B_r(m, w)|. \tag{6.6}$$

Secondly, partition the set $\mathscr{A} \cap B(m)$ into two subsets $\mathscr{A} \cap B_l(m, w)$ and $\mathscr{A} \cap B_r(m, w)$:
$$\mathscr{A} \cap B(m) = [\mathscr{A} \cap B_l(m, w)] \cup [\mathscr{A} \cap B_r(m, w)].$$
Due to $[\mathscr{A} \cap B_l(m, w)] \cap [\mathscr{A} \cap B_r(m, w)] = \varnothing$, the following partition of $|\mathscr{A} \cap B(m)|$ also holds:
$$|\mathscr{A} \cap B(m)| = |\mathscr{A} \cap B_l(m, w)| + |\mathscr{A} \cap B_r(m, w)|. \tag{6.7}$$

Now deduce the comparing result of both equations (3.1) and (6.5).

**Theorem 6.4.** *Suppose that $N$ is a large even number, $w < 0.5 N$ and $m \ge 2$ are real numbers. The following estimate holds*:
$$|\mathscr{A} \cap B_l(m, w)| = 2\alpha(m, N) \frac{c(N)}{\ln N} |B_l(m, w)| \left(1 + O\left(\frac{1}{\ln^{1/14} N}\right)\right). \tag{6.8}$$

*Proof.* From equations (3.1), (6.5) and (6.7) one obtains
$$|\mathscr{A} \cap B_l(m, w)| = |\mathscr{A} \cap B(m)| - |\mathscr{A} \cap B_r(m, w)|$$
$$= 2\alpha(m, N) \frac{c(N)}{\ln N} |B(m)| \left(1 + O\left(\frac{1}{\ln^{1/14} N}\right)\right)$$
$$- 2\alpha_1(m, N) \frac{c(N)}{\ln N} |B_r(m, w)| \left(1 + O\left(\frac{1}{\ln^{1/14} N}\right)\right). \tag{6.9}$$

From equation (6.6), equation (6.9) can be written in

$$|\mathcal{A} \cap B_l(m,w)| = 2\alpha(m,N)\frac{c(N)}{\ln N}|B_l(m,w)|\left(1+O\left(\frac{1}{\ln^{1/14} N}\right)\right)+$$

$$+ 2\alpha(m,N)\frac{c(N)}{\ln N}|B_r(m,w)|\left(1+O\left(\frac{1}{\ln^{1/14} N}\right)\right)-$$

$$- 2\alpha_1(m,N)\frac{c(N)}{\ln N}|B_r(m,w)|\left(1+O\left(\frac{1}{\ln^{1/14} N}\right)\right). \quad (6.10)$$

Note that

$$B_r(m,w) \cap [\mathcal{A} \cap B_l(m,w)] = \varnothing,$$

$$B_r(m,w) \cap B_l(m,w) = \varnothing,$$

therefore, $|\mathcal{A} \cap B_l(m,w)|$ and $|B_l(m,w)|$ are independent from $|B_r(m,w)|$, and then there must be

$$2\alpha(m,N)\frac{c(N)}{\ln N}|B_r(m,w)|\left(1+O\left(\frac{1}{\ln^{1/14} N}\right)\right)=$$

$$= 2\alpha_1(m,N)\frac{c(N)}{\ln N}|B_r(m,w)|\left(1+O\left(\frac{1}{\ln^{1/14} N}\right)\right)$$

in the equation (6.10), thus Theorem 6.4 holds.

The proof of the Proposition {1, 1} is then complete:

**Theorem 6.5.** *Let $D(m,N)$ be the number of the primes not exceeding $N^{2/m}$ in the set $\mathcal{A}$. The following estimate holds when $m \geq 2$:*

$$D(m,N) = m\alpha(m,N)c(N)\frac{N^{2/m}}{\ln^2 N}\left(1+O\left(\frac{1}{\ln^{1/14} N}\right)\right). \quad (6.11)$$

*where $\alpha(m,N)$ is a real number with $\frac{f(m/2)}{F(m)} \leq \alpha(m,N) \leq \frac{F(m/2)}{f(m)}$ and $\lim_{m\to\infty}\alpha(m,N)=1$.*

*Proof.* Put $w = N^{2/m}$ in the set $B_l(m,w)$ of equation (6.8), noting that

$$|B_l(m,N^{2/m})| = \pi(N^{2/m}) - \pi(N^{1/m}) = \frac{mN^{2/m}}{2\ln N}\left(1+O\left(\frac{1}{\ln N}\right)\right)$$

because the integers in the set $B_l(m,N^{2/m})$ are all primes; from equation (6.8) one obtains

$$D(m,N) = |\mathcal{A} \cap B_l(m,N^{2/m})| + O(\pi(N^{1/m}))$$

$$= 2\alpha(m,N)\frac{c(N)}{\ln N}|B_l(m,N^{2/m})|\left(1+O\left(\frac{1}{\ln^{1/14} N}\right)\right)$$

$$= m\alpha(m, N)c(N)\frac{N^{2/m}}{\ln^2 N}\left(1+O\left(\frac{1}{\ln^{1/14} N}\right)\right).$$

The theorem holds.

The significance of Theorem 6.5 lies in the fact that Proposition $\{1, 1\}$ can be proved for any $m \geq 5$. For example, $\alpha(5, N) \geq \frac{f(5/2)}{F(5)} = 0.5767$ when $m = 5$, it follows, from Theorem 6.5, that

$$\left|\mathscr{A} \cap \pi\right| \geq D(5, N) \geq 2.8835 \frac{c(N)N^{2/5}}{\ln^2 N}\left(1+O\left(\frac{1}{\ln^{1/14} N}\right)\right) > 0.$$

Thus Proposition $\{1, 1\}$ is true.

The uncertainty in the above results produced by the unknown value of the coefficient $\alpha(m, N)$ can be eliminated by proving that $\alpha(m, N) = 1$ for any $m \geq 2$, see below.

**Theorem 6.6.** *Suppose that $N$ is a large even number. The following equation holds for any real number $m \geq 2$:*

$$\alpha(m, N) = 1. \tag{6.12}$$

*Proof.* When $m \geq 2$, Theorem 6.5 gives

$$D(m, N) = m\alpha(m, N)c(N)\frac{N^{2/m}}{\ln^2 N}\left(1+O\left(\frac{1}{\ln^{1/14} N}\right)\right). \tag{6.13}$$

On the other hand, when $m' \geq 2$ with $m/2 < m' < m$, Theorem 6.4 gives

$$\left|\mathscr{A} \cap B_l(m', w)\right| = 2\alpha(m', N)\frac{c(N)}{\ln N}\left|B_l(m', w)\right|\left(1+O\left(\frac{1}{\ln^{1/14} N}\right)\right). \tag{6.14}$$

Put $w = N^{2/m}$ in the sets $\mathscr{A} \cap B_l(m', w)$ and $B_l(m', w)$ of equation (6.14), and calculating the first approximation of $D(m, N)$ from equation (6.14). Noting that $\frac{2}{m'} > \frac{2}{m} > \frac{1}{m'}$ holds when $m' \geq 2$ with $m/2 < m' < m$, which implies that the integers in the sets $\mathscr{A} \cap B_l(m', N^{2/m})$ and $B_l(m', N^{2/m})$ are all primes, thereout one obtains

$$\left|B_l(m', N^{2/m})\right| = \pi(N^{2/m}) - \pi(N^{1/m'})$$

$$= \frac{mN^{2/m}}{2\ln N}\left(1+O\left(\frac{1}{\ln N}\right)\right),$$

then, from Theorem 6.5, one obtains

$$D(m, N) = \left|\mathscr{A} \cap B_l(m', N^{2/m})\right| + O(\pi(N^{1/m'}))$$

$$= 2\alpha(m', N)\frac{c(N)}{\ln N}\left|B_l(m', N^{2/m})\right|\left(1 + O\left(\frac{1}{\ln^{1/14} N}\right)\right)$$

$$= m\alpha(m', N)\frac{c(N)N^{2/m}}{\ln^2 N}\left(1 + O\left(\frac{1}{\ln^{1/14} N}\right)\right). \qquad (6.15)$$

It is clear, by comparing the $D(m, N)$ in equation (6.15) with that in equation (6.13), that $\alpha(m', N) = \alpha(m, N)$ since $D(m, N)$ and its first approximation are the same in both cases. The fact that $\alpha(m', N) = \alpha(m, N)$ must hold for any even integer $N$ and any real number pair $m \geq 2$ and $m' \geq 2$ with $m/2 < m' < m$ implies that the value of $\alpha(m, N)$ is definitely not related to $m$. It is well known that $\lim_{m \to \infty} \alpha(m, N) = 1$, and so we must have $\alpha(m, N) = 1$ for any real number $m \geq 2$. Therefore, the theorem holds.

The first approximation to the Goldbach problem can be obtained as follows:

**Theorem 6.7.** *Suppose that $N$ is a large even number. The following estimate holds*:

$$\left|\mathscr{A} \cap \pi\right| = 2c(N)\frac{N}{\ln^2 N}\left(1 + O\left(\frac{1}{\ln^{1/14} N}\right)\right), \qquad (6.16)$$

*where $\mathscr{A} \cap \pi$ is the number of representations of $N$ as the sum of two primes.*

*Proof.* When $m = 2$ and $\alpha(m, N) = 1$, Theorem 3.1 gives

$$\left|\mathscr{A} \cap B(2)\right| = 2\frac{c(N)}{\ln N}\left|B(2)\right|\left(1 + O\left(\frac{1}{\ln^{1/14} N}\right)\right).$$

Noting that

$$\left|\mathscr{A} \cap \pi\right| = \left|\mathscr{A} \cap B(2)\right| + O(\pi(N^{1/2})),$$

$$\left|B(2)\right| = \pi(N) - \pi(N^{1/2}) = \frac{N}{\ln N}\left(1 + O\left(\frac{1}{\ln N}\right)\right),$$

Clearly, the theorem holds.

Theorem 6.7 is precisely the answer to the Goldbach problem! We see that the problem can be solved by the sieve method which only relate to a few technical skill.

**7. The problem of twin primes**

The above results can equally well be applied to the problem of twin primes.

Suppose that $b$ is a positive even integer and $x$ is a real number. Introduce the set $\mathscr{C}_b(x)$:

$$\mathscr{C}_b(x) = \{t|\ b < t \leq x,\ t \in \pi,\ t - b \in \pi\}.$$

Then one has

**Theorem 7.1.** $$\left|\mathscr{C}_b(x)\right| \sim 2c(b)\frac{x}{\ln^2 x}. \tag{7.1}$$

The above results for the Goldbach problem and the problem of twin primes are entirely consistent with the conjectures of Hardy-Littlewood [19].

**Acknowledgments**
The author is sincerely grateful to Dr. Ji-guo Chu, Dr. Francis R. Austin, and Prof. Ke-ren Liao for their zealous help and support.